\documentclass{article}

\usepackage{makeidx}  
\usepackage{graphicx} 
\usepackage{subeqnar} 
\usepackage{multicol} 
\usepackage{cropmark} 
\makeindex            

\usepackage{amsmath,amssymb}   
\usepackage{amsfonts,fancybox}

\newcommand{\real}{{\mathbb{R}}}

\newcommand{\symprod}[2]{\langle #1: #2\rangle}
\newcommand{\symprodmed}[2]{\Big\langle #1 : #2\Big\rangle}
\newcommand{\pder}[2]{\frac{\partial #1}{\partial #2}}
\newcommand{\eqdef}{\triangleq}
\newcommand{\mcI}{{\mathcal{I}}}
\newcommand{\mcU}{{\mathcal{U}}}
\newcommand{\overl}{\overline}
\newcommand{\eps}{\epsilon}

\newcommand{\spn}{\operatorname{span}}
\newcommand{\lift}{^{\text{lift}}}

\newtheorem{theorem}{Theorem}[section]
\newtheorem{remark}[theorem]{Remark}

\begin{document}
\title{\vspace{-5ex}\hfill        \texttt{\begin{minipage}{.7\textwidth}\small
Worskshop on Control Problems in Robotics and Automation, Las Vegas, December
14,  2002\end{minipage}}\\ [3ex]  Motion  planning and  control problems  for
underactuated robots}

\author{Sonia Mart{\'\i}nez\\
     Escuela Universitaria Polit\'ecnica\\
     Universidad Polit\'ecnica de Catalu\~na\\
     Av. V. Balaguer s/n\\
     Vilanova i la Geltr\'u 08800, Spain
\and Jorge Cort\'es, Francesco Bullo\\
     Coordinated Science Laboratory\\
     University of Illinois at Urbana-Champaign\\
     1308 W. Main St\\
     Urbana, IL 61801, USA}

\markboth{S. Mart{\'\i}nez, J. Cort\'es, F. Bullo}{A ver}

\maketitle              

\begin{abstract}
  Motion planning and control are key problems in a collection of robotic
  applications including the design of autonomous agile vehicles and of
  minimalist manipulators.  These problems can be accurately formalized
  within the language of affine connections and of geometric control theory.
  In this paper we overview recent results on kinematic controllability and
  on oscillatory controls. Furthermore, we discuss theoretical and practical
  open problems as well as we suggest control theoretical approaches to them.
\end{abstract}

\section{Motivating problems from a variety of robotic applications} 

The research in Robotics is  continuously exploring the design of novel, more
reliable and agile  systems that can provide more  efficient tools in current
applications  such  as factory  automation  systems,  material handling,  and
autonomous robotic applications, and  can make possible their progressive use
in areas such as medical and social assistance applications.

Mobile  Robotics,  primarily  motivated   by  the  development  of  tasks  in
unreachable  environments, is  giving way  to new  generations  of autonomous
robots in  its search for new  and ``better adapted''  systems of locomotion.
For  example, traditional  wheeled  platforms have  evolved into  articulated
devices  endowed with  various types  of wheels  and suspension  systems that
maximize their traction and the robot's ability to move over rough terrain or
even climb  obstacles.  The types of  wheels that are  being employed include
passive  and powered  castors,  ball-wheels or  omni-directional wheels  that
allow  a  high  accuracy  in  positioning and  yet  retain  the  versatility,
flexibility and  other properties  of wheels.  A  rich and  active literature
includes                  (i)                 various                 vehicle
designs~\cite{SO-JA-RS:00,FGP-SMK:94,SS-JA-JD:95,MW-HA:97},      (ii)     the
automated    guided    vehicle    ``OmniMate''~\cite{JBo:00},    (iii)    the
roller-walker~\cite{GE-SH:99}  and other dexterous  systems~\cite{SH:00} that
change  their internal  shape and  constraints  in response  to the  required
motion sequence, and (iv) the omni-directional platform in~\cite{RH-OK:00}.

Other types of remotely controlled autonomous vehicles that are increasingly
being employed in space, air and underwater applications include
submersibles, blimps, helicopters, and other crafts.  More often than not
they rely on innovative ideas to affect their motion instead of on classic
design ideas. For example, in underwater vehicle applications, innovative
propulsion systems such as shape changes, internal masses, and momentum
wheels are being investigated.  Fault tolerance, agility, and maneuverability
in low velocity regimes, as in the previous example systems, are some of the
desired capabilities.

\begin{figure}
  \includegraphics[width=.35\textwidth]{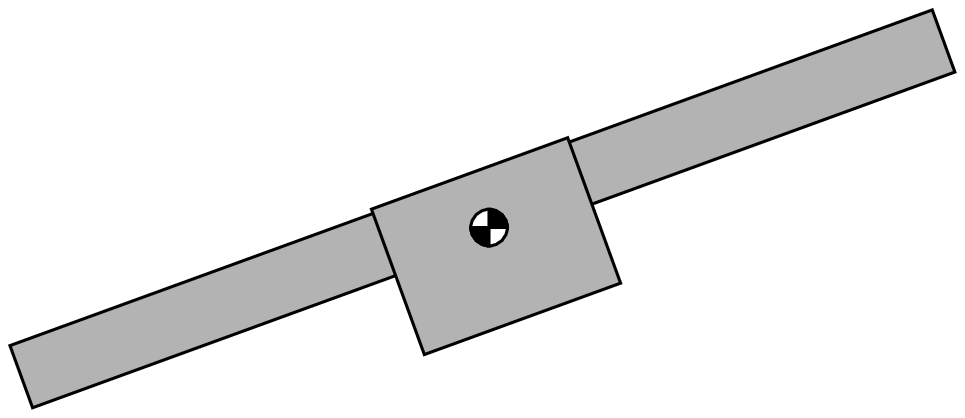}
  \includegraphics[width=.25\textwidth]{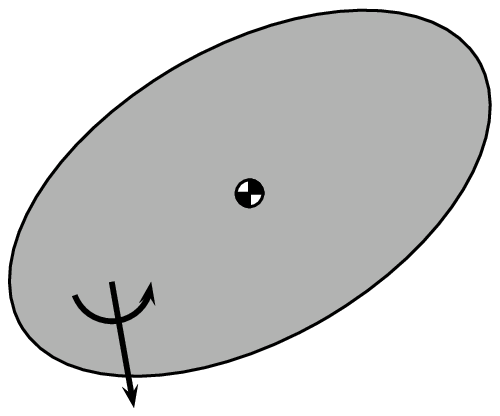}
  \includegraphics[width=.35\textwidth]{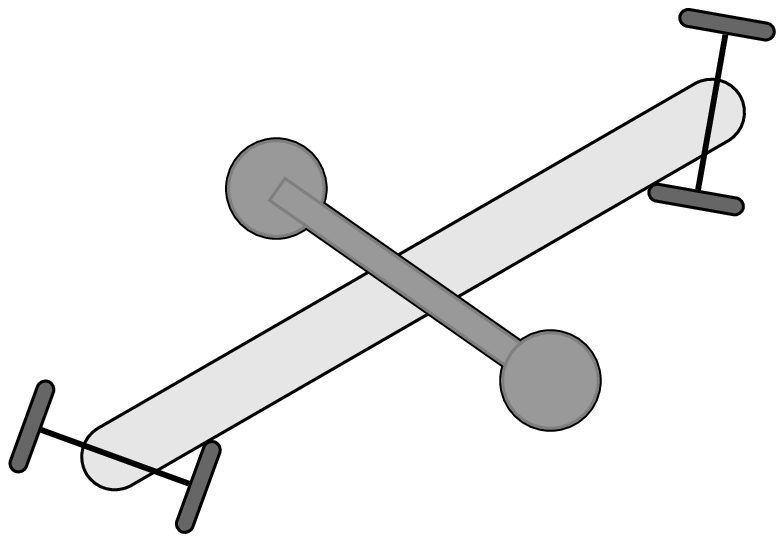}
  \caption{Underactuated robots appear in a variety of environments. From left
    to right, a planar vertical take-off and landing (PVTOL) aircraft model,
    a horizontal model of a blimp and the snakeboard.}
  \label{fig:examples}
\end{figure}

A growing field in Mobile Robotics is that of \emph{biomimetics}. The idea of
this  approach is  to obtain  some of  the robustness  and  adaptability that
biological systems have refined  through evolution. In particular, biomimetic
locomotion  studies  the  periodic  movement patterns  or  \emph{gaits}  that
biological systems undergo  during locomotion and then takes  it as reference
for the design  of the mechanical counterpart. In other  cases, the design of
to  robots without physical  counterpart is  inspired by  similar principles.
Robotic locomotion systems include  the classic \emph{bipeds and multi-legged
robots}  as well  as swimming  snake-like  robots and  flying robots.   These
systems find potential applications  in harsh or hazardous environments, such
as  under  deep or  shallow  water, on  rough  terrain  (with stairs),  along
vertical  walls or  pipes  and  other environments  difficult  to access  for
wheeled robots.  Specific examples  in the literature include hyper-redundant
robots~\cite{GC-JWB:94,SH:93},                                             the
snakeboard~\cite{ADL-JPO-RMM-JWB:94,JPO-JWB:98},  the  $G$-snakes and  roller
racer          models          in~\cite{PSK-DPT:94,PSK-DPT:01},          fish
robots~\cite{NK-TI:98,SDK-RJM-CTA-RMM-JWB:98},                             eel
robots~\cite{JCJ-SK-JA:95,KAM-JPO:99a},     and    passive     and    hopping
robots~\cite{JKH-MHR:90,TM:90,MHR:86}.

All this set of emerging robotic applications have special characteristics 
that pose new challenges in motion planning. Among them, we highlight:

\paragraph{Underactuation.} This could be owned to a design choice: nowadays
low weight and fewer actuators must perform the task of former more expensive
systems.  For example, consider a manufacturing environment where robotic
devices perform material handling and manipulation tasks: automatic planning
algorithms might be able to cope with failures without interrupting the
manufacturing process.  Another reason why these systems are underactuated is
because of an unavoidable limited control authority: in some locomotion
systems it is not possible to actuate all the directions of motion.  For
example, consider a robot operating in a hazardous or remote environment
(e.g., aerospace or underwater), an important concern is its ability to
operate faced with a component failure, since retrieval or repair is not
always possible.

\paragraph{Complex dynamics.} In these control systems, the drift 
plays a key role.  Dynamic effects must necessarily be taken into account,
since kinematic models are no longer available in a wide range of current
applications. Examples include lift and drag effects in underwater vehicles,
the generation of momentum by means of the coupling of internal shape changes
with the environment in the eel robot and the snakeboard, the dynamic
stability properties of walking machines and nonholonomic wheeled platforms,
etc.

\paragraph{Current limitations of motion algorithms.}
Most of the work on motion planning has relied on assumptions that are no
longer valid in the present applications.  For example, one of these is that
(wheeled) robots are kinematic systems and, therefore, controlled by velocity
inputs. This type of models allows one to design a control to reach a desired
point and then immediately stop by setting the inputs to zero. This is
obviously not the case when dealing with complex dynamic models.

Another  common assumption  is  the one  of  fully actuation  that allows  to
decouple  the  motion  planning  problem  into  path  planning  (computational
geometry)  and then  tracking.  For  underactuated  systems, this  may be  not
possible because we  may be obtaining motions in the  path planning stage that
the  system can  not  perform in  the  tracking step  because  of its  dynamic
limitations.

\begin{figure}[tbh]
    \includegraphics[width=.42\textwidth]{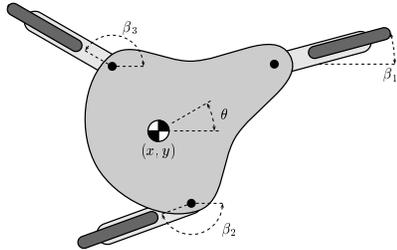}
    \caption{Vertical  view of  an  omni-directional robotic  platform with  6
      degrees of freedom and 3 nonholonomic
      constraints~\cite{GC-GB-BDAN:96,RH-OK:00}.  This device is capable of
      highly accurate positioning, high payloads, and high speed motion.  In
      its fully actuated configuration, the robot is endowed with 6 motors at
      the three wheels and at the three joints $(\beta_1,\beta_2,\beta_3)$.
      However, underactuated configurations can arise because of failures or
      intentional design.}
\label{fig:omni}
\end{figure}

Furthermore, motion planning and optimization problems for these systems are
nonlinear, non-convex problems with exponential complexity in the dimension
of the model. These issues have become increasingly important due to the high
dimensionality of many current mechanical systems, including flexible
structures, compliant manipulators and multibody systems undergoing
reconfiguration in space.

\paragraph{Benefits that would result from better motion planning algorithms
for underactuated systems.}  From a practical perspective, there  are at least
two advantages to designing controllers for underactuated robotic manipulators
and vehicles. First, a fully actuated system requires more control inputs than
an underactuated  system, which means  there will have  to be more  devices to
generate the necessary forces.  The  additional controlling devices add to the
cost and  weight of  the system.   Finding a way  to control  an underactuated
version  of the system  would improve  the overall  performance or  reduce the
cost. The second practical reason  for studying underactuated vehicles is that
underactuation  provides  a backup  control  technique  for  a fully  actuated
system.   If a  fully  actuated system  is  damaged and  a  controller for  an
underactuated system is  available, then we may be  able to recover gracefully
from  the failure.   The underactuated  controller may  be able  to  salvage a
system that would otherwise be uncontrollable. 


\section{Mathematical unifying approach to the modeling of robotic 
systems}\label{sec:modeling}

Most of the robotic devices we have mentioned so far can be characterized by
their special Lagrangian structure. They usually exhibit symmetries and
their motion is constrained by the environment where they operate. In the
following, we introduce a general modeling language for underactuated robotic
systems.

Let $q=(q^1,\ldots,q^n)\in  Q$ be the  configuration of the  mechanical system
and consider the control equations:
\begin{gather}
  \label{eq:mechsys}
  \ddot{q}^i  + \Gamma_{jk}^i(q)\dot{q}^j\dot{q}^k  = -  M^{ij}\pder{V}{q^j}
  + k^i_j(q) \dot{q}^j + Y_1^i(q) u_1 + \ldots + Y_m^i(q) u_m \, ,
\end{gather}
where the summation convention is in place for the indices $j,k$ that run from
$1$ to $n$, and
\begin{enumerate}
\item $V:Q \rightarrow \real$ corresponds to potential energy, and $k^i_j(q) \dot{q}^j$ corresponds to damping forces,
\item  $\{\Gamma_{jk}^i: i,j,k=1,\ldots,n\}$  are  $n^3$ Christoffel  symbols,
derived from $M(q)$, the inertia matrix defining the kinetic energy, according
to
  \begin{equation*}
    \Gamma_{{i}{j}}^k = \frac{1}{2} M^{mk}\left(
      \pder{M_{m{j}}}{q^{i}} +
      \pder{M_{m{i}}}{q^{j}}-\pder{M_{{i}{j}}}{q^m} \right) \, ,          
  \end{equation*}
  where  $M^{mk}$ is the $(m,k)$ component of $M^{-1}$, and,
\item $\{F_a: a=1,\ldots,m\}$ are the $m$ input co-vector fields, and $\{Y_a =
  M^{-1} F_a: a=1,\ldots,m\}$ are the $m$ input vector fields.
\end{enumerate}

Underactuated systems have fewer control actuators, $m$, than degrees of
freedom $n>m$.  Other limitations on the control signals $u_a$ might be
present, e.g., actuators might have magnitude and rate limits, or they might
only generate unilateral or binary signals (e.g., thrusters in satellites).

The notion of affine connection provides a coordinate-free means of
describing the dynamics of robotic systems.  Given two vector fields $X,Y$,
the \emph{covariant derivative} of $Y$ with respect to $X$ is the third
vector field $\nabla_XY$ defined via
\begin{equation}\label{cov:der:Christoffels}
  (\nabla_X Y)^i  = \frac{\partial Y^i}{\partial q^j} X^j  + \Gamma^i_{jk} X^j
  Y^k .
\end{equation}
The operator $\nabla$ is called the \emph{affine connection} for the
mechanical system in equation~\eqref{eq:mechsys}.  We write the
Euler-Lagrange equations for a system subject to a conservative force $Y_0$,
a damping force $k(q)\dot{q}$ and $m$ input forces as:
\begin{equation}
  \label{eq:mechsys:affine}
 \nabla_{\dot{q}}\dot{q}  = Y_0(q)  + k(q)(\dot{q})  + \sum_{a=1}^{m}  Y_a (q)
 u_a(t).
\end{equation}
Equation~\eqref{eq:mechsys:affine} is a coordinate-free version of
equation~\eqref{eq:mechsys}. A crucial observation is the fact that systems
subject nonholonomic constraints can also be modeled by means of affine
connections. In the interest of brevity, we refer to~\cite{FB-MZ:01d,ADL:97a}
for the exposition of this result and the explicit expression of the
Christoffel symbols corresponding to the Lagrange-d'Alembert equations.

\subsubsection*{The       homogeneous       structure      of       mechanical
  systems.}\label{subsec:Lie}

The     fundamental     structure     of     the    control     system     in
equation~\eqref{eq:mechsys:affine}  is  the   polynomial  dependence  of  the
various  vector fields on  the velocity  variable $\dot{q}$.   This structure
affects the Lie bracket computations involving input and drift vector fields.
The system~\eqref{eq:mechsys:affine}  is written in  first order differential
equation form as
\begin{eqnarray*}
  \frac{d}{dt} \begin{bmatrix} q \\ \dot{q}  \end{bmatrix} =
  \begin{bmatrix} \dot{q} \\ -\Gamma(q,\dot{q}) +Y_0(q) 
    +k(q)(\dot{q})\end{bmatrix}
  + \sum_{a=1}^m \begin{bmatrix}0\\Y_a\end{bmatrix}u_a(t)
\end{eqnarray*}
where $\Gamma(q,\dot{q})$ is the vector with $i$th component
$\Gamma^i_{jk}(q)\dot{q}^j\dot{q}^k$. Also, if  $x=(q,\dot{q})$,
\begin{eqnarray*}
Z(x) = \begin{bmatrix} \dot{q} \\ -\Gamma(q,\dot{q}) \end{bmatrix} ,
\quad
 Y_a\lift(x) \eqdef  \begin{bmatrix} 0 \\  Y_a(q) \end{bmatrix} ,
\quad \text{and} \quad
 k\lift(x) \eqdef  \begin{bmatrix} 0 \\  k(q)(\dot{q}) \end{bmatrix} ,
\end{eqnarray*}
the control system is rewritten as
\begin{equation*}
\dot{x} = Z(x) + Y_0\lift(x) + k\lift(x)  + \sum_{a=1}^m Y_a\lift(x)
u_a(t) \, .
\end{equation*}

Let $h_i(q,\dot{q})$ be the set of scalar functions on $\real^{2n}$ which are
arbitrary    functions   of    $q$\/   and    homogeneous    polynomials   in
$\{\dot{q}^1,\dots,\dot{q}^n\}$ of degree $i$.  Let $\mathcal{P}_i$ \/ be the
set of  vector fields  on $\real^{2n}$ whose  first $n$ components  belong to
$h_i$  and whose second  $n$ components  belong to  $h_{i+1}$.  We  note that
these   notions    can   also   be    defined   on   a    general   manifold,
see~\cite{FB-ADL:00d}.

We are now ready to present two simple ideas. First, all the previous vector
fields are homogeneous polynomial vector fields for some specific value of
$i$.  Indeed, $Z \in \mathcal{P}_1$, $k\lift \in \mathcal{P}_{0}$, and
$Y_a\lift \in \mathcal{P}_{-1}$.  Second, since the Lie bracket between a
vector field in~$\mathcal{P}_i$ and a vector field in~$\mathcal{P}_j$ belongs
to~$\mathcal{P}_{i+j}$, any Lie bracket of the given relevant vector fields
remains a homogeneous polynomial.  In other words, the set of homogeneous
vector fields is closed under the operation of Lie bracket.

A consequence of this analysis is the definition of symmetric product of
vector fields. We define the \emph{symmetric product} between $Y_b$ and $Y_a$
as the vector field $\symprod{Y_a}{Y_b} = \symprod{Y_b}{Y_a}$ given by
\begin{align*}
  \symprod{Y_b}{Y_a}^i  =   \symprod{Y_a}{Y_b}^i
  &= \frac{\partial Y_a^i}{\partial q^j}Y_b^j
  +\frac{\partial Y_b^i}{\partial q^j}Y_a^j + \Gamma^i_{jk}
  \left(Y_a^jY_b^k + Y_a^kY_b^j\right).
\end{align*}
Straightforward  computations show  that  $\symprod{Y_a}{Y_b}^{\text{lift}} =
[Y_b^{\text{lift}},[Z_g,Y_a^{\text{lift}}]]$. This operation plays a key role
in nearly  all the  control problems associated  with this class  of systems:
nonlinear       controllability~\cite{JC-SM-FB:01m,ADL-RMM:95c},      optimal
control~\cite{MC-FSL-PEC:95,ADL:99a},             dynamic            feedback
linearization~\cite{MR-RMM:96b},   algorithms   for   motion   planning   and
stabilization~\cite{FB-NEL-ADL:00,SM-JC:02,JPO:00a}, etc.

\subsubsection*{A series expansion for the forced evolution starting from
  rest.} \label{sec:series-expansion}

The     homogeneous      structure     of     the      mechanical     control
system~\eqref{eq:mechsys:affine},  together with  the symmetric  product, set
the  basis to establish  the following  description of  the evolution  of the
system       trajectories        starting       with       zero       initial
velocity~\cite{FB:99a,JC-SM-FB:01m}.  Assume  no potential or  damping forces
are present in the system. Let $Y(q,t) = \sum_{a=1}^m Y_a(q) u_a(t)$.  Define
recursively the vector fields $V_k$ by
\begin{align*} 
  V_1(q,t) &= \int_0^tY(q,s)ds  \, , \; \; 
  V_k(q,t) = - \frac{1}{2} \sum_{j=1}^{k-1} \;
  \int_0^t\!\symprodmed{V_j(q,s)}{\null\; V_{k-j}(q,s)} ds.
\end{align*}  
Then, the solution $q(t)$ of equation~\eqref{eq:mechsys:affine} satisfies
\begin{equation} \label{eq:thm3}
  \dot{q}(t) = \sum_{k=1}^{+\infty} {V}_k(q(t),t),
\end{equation}
where the series converges absolutely and uniformly in a neighborhood
of $q_0$ and over a fixed time interval $t\in [0,T]$. This series expansion
provides a means of describing the open-loop response of the system to any
specific forcing. As we will see below, it plays a key role in several motion
planning and control strategies for underactuated robots.


\section{Existing results on planning for underactuated systems}

To design planning algorithms  for underactuated robotic systems, we advocate
an  integrated approach  based  on modeling,  system design,  controllability
analysis,  dexterity,  manipulability,  and  singularities.   These  analysis
concepts are  fundamental for robust  planning algorithms that do  not solely
rely  on   randomization  or  nonlinear  programming.   We   do  not  suggest
closed-form planning algorithms, rather  we envision methods that combine the
best features of formal analysis and of numerical algorithms.

For reasons of space, we cannot present a detailed account of all existing
results on motion planning for underactuated systems, and not even of the
results obtained within the modeling approach proposed in
Section~\ref{sec:modeling}.  Therefore, we focus on two specific control
methodologies for motion planning: decoupled planning algorithms for
kinematically controllable systems, and approximate inversion algorithms based
on oscillatory controls.

Section~\ref{subsec:kinematic-controllability} reviews \emph{decoupled
  planning algorithms} that exploit certain differential geometric properties
to reduce the complexity of the motion planning problem (still to be solved
via numerical algorithms).  The notion of \emph{kinematic controllability} is
extremely effective: trajectory planning decouples from being a problem on
a~$2n$\ dimensional space to an~$n$\ dimensional space. Furthermore, various
state constraints can be neglected in the reduced space.  For systems that
are not kinematically controllable and that require oscillatory controls to
locomote, Section~\ref{subsec:approximate-inversion} presents motion planning
algorithms based on \emph{approximate inversion}.  Both design methods are
closely related to recent results on nonlinear
controllability~\cite{FB-KML:01a,ADL-RMM:95c}, power series
expansions~\cite{FB:99a,JC-SM-FB:01m}, two time-scales coordinate-free
averaging~\cite{FB:99b,SM-JC-FB:01e}, and nonlinear inversion
algorithms~\cite{FB-NEL-ADL:00,SM-JC:02}.

The strengths of this methodology  are as follows. Both methodologies provide
solutions to the corresponding problems,  i.e., point to point and trajectory
planning.  These analytic results do not rely on non-generic assumptions such
as  feedback   linearization,  nilpotency  or  flatness.    The  results  are
coordinate-free and hence widely applicable, e.g., to aerospace or underwater
robotics settings.   Both methods are  \emph{consistent}, \emph{complete} and
constructive (consistent  planners recover the known  solutions available for
linear and nilpotent systems, and  complete planners are guaranteed to find a
local solution for any nonlinearly controllable system).

\subsection{Kinematic controllability for underactuated robots}
\label{subsec:kinematic-controllability}
The following decoupling methodology was proposed in~\cite{FB-KML:01a} to
reduce the complexity of the motion planning problem.  The method is
constructive (only quadratic equations and no PDEs are involved) and
physically intuitive.

\begin{figure}
  \includegraphics[width=.4\textwidth]{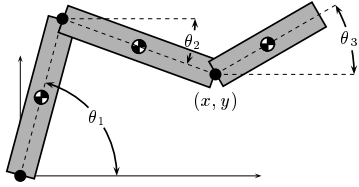}
  \caption{A three-revolute-joints device. It can be proven \cite{FB-KML:01a}
  that  any  two-actuator  configuration  of  this  system  is  kinematically
  controllable, i.e., one can always  find two decoupling vector fields whose
  involutive closure is full-rank.}
  \label{fig:three-links}
\end{figure}

We  consider as  a  motivating example  a  common pick-\&-place  manipulator:
Fig.~\ref{fig:three-links} shows  a vertical view  of a three-revolute-joints
device.   We investigate planning  schemes for  this system  when one  of its
three motors  is either failed or  missing.  We present a  decoupling idea to
reduce  the complexity  of the  problem:  instead of  searching for  feasible
trajectories of  a dynamic  system in~$\real^6$, we  show how it  suffices to
search  for  paths  of  a  simpler, kinematic  (i.e.,  driftless)  system  in
$\real^3$.

A   curve   $\gamma:[0,T]\mapsto   Q$\    is   a   controlled   solution   to
equation~\eqref{eq:mechsys}  if there  exist inputs  $u_a:[0,T]\to\real$\ for
which  $\gamma$ solves~\eqref{eq:mechsys}.  To  avoid the  difficult task  of
characterizing all controlled  solutions of the system~\eqref{eq:mechsys}, we
focus on curves satisfying $\dot{\gamma}=\dot{s}(t)X(\gamma)$, where $X$ is a
vector  field   on  $Q$,  and  where   the  map  $s:[0,T]\to   [0,1]$\  is  a
``time-scaling''  parameterization   of~$\gamma$.  Such  curves   are  called
\emph{kinematic  motions}.

We call $V$\ a \emph{decoupling vector field} if all curves $\gamma$\ 
satisfying $\dot{\gamma}=\dot{s}(t)V(\gamma)$\ for any time scaling $s$, are
kinematic motions.  This definition is useful for three reasons. First, $V$\ 
is decoupling if and only if $V$\ and $\nabla_VV$\ are linear combinations in
$\{Y_1,\dots,Y_m\}$.  Second, decoupling vector fields can be computed by
solving $(n-m)$ quadratic equations.  Third, if enough decoupling vector
fields, say $V_1,\dots,V_p$, are available to satisfy the LARC, we call the
system \emph{kinematically controllable}.  In the latter case, we can plan
motions for the kinematic system $\dot{q} = \sum_{a=1}^p w_a(t) V_a(q)$, and
they will automatically be controlled curves for the original
system~\eqref{eq:mechsys}.

\subsection{Approximate inversion via small amplitude and oscillatory controls}
\label{subsec:approximate-inversion}

As in  the previous section, the  objective is to design  motion planning and
stabilization schemes for underactuated systems.  We propose perturbation and
inversion methods as widely applicable approaches to solve point to point and
trajectory  planning  problems.  Let  us   regard  the  flow  map  $\Phi$  of
equation~\eqref{eq:mechsys:affine} over a finite time interval as a map from
the  input functions $u_i:[0,T]\to\real$\  to the  target state  $x(T)$.  The
ideal algorithm for point-to-point planning computes an exact (right) inverse
$\Phi^{-1}$   of~$\Phi$.   Unfortunately,   closed   form   expressions   for
$\Phi^{-1}$\  are  available  only assuming  \emph{non-generic}  differential
geometric  conditions (e.g., the  system needs  to be  feedback linearizable,
differentially  flat,   or  nilpotent).   Instead  of   aiming  at  ``exact''
solutions,  we focus  on computing  an \emph{approximate  inverse  map} using
perturbation methods  such as power  series expansions and  averaging theory.
Although  these tools  are  only approximate,  the  resulting algorithms  are
consistent and complete.
 
\subsubsection*{Oscillatory (high frequency, high amplitude) controls for
  trajectory planning.}  

We present the approach in three steps and refer to~\cite{SM-JC-FB:01e} for
all the details.  As first step, we present a recent coordinate-free
averaging result.  Let $0<\eps\ll 1$. Assume the control inputs are of the
form
\begin{equation*}
  u_i = \frac{1}{\eps}  u_i\left(\frac{t}{\eps},t\right),
\end{equation*}
and assume they are $T$-periodic and zero-mean in the first variable.  Define
the averaged multinomial iterated integrals of $u_1,\dots,u_m$ as
\begin{equation*}
  \overl{\mcU}_{k_1,  \dots,  k_m}(t)   =  \frac{T^{-1}}{k_1!  \dots  k_{m}!}
  \int_0^T\!\!      \left(\int_0^s\!\!u_1(\tau,t)      d\tau     \right)^{k_1}
  \hspace*{-5pt}  \dots  \left(\int_0^s\!\!  u_m(\tau,t)  d\tau  \right)^{k_m}
  \hspace*{-5pt} ds \, .
\end{equation*}
Let $a,b,c$ take value in $\{1,\dots,m\}$.  Let $\vec{k}_a$ (resp.
$\vec{k}_{ab}$) denote the tuple $(k_1,\dots,k_m)$ with $k_c = \delta_{ca}$
(resp.  $k_c =\delta_{ca} + \delta_{cb}$). Then, over a finite time $q(t) =
r(t) + O(\eps)$, as $\eps \rightarrow 0$, where $r(t)$ satisfies
\begin{align}\label{eq:mechsys:averaged}
&  \nabla_{\dot{r}}\dot{r} =  Y_{0}(r) +  k(r)(\dot{r}) +  \sum_{a=1}^m \left(
\frac{1}{2} \overl{\mcU}_{\vec{k}_a}^2 (t) - \overl{\mcU}_{\vec{k}_{aa}}(t)\right)
\symprod{Y_a}{Y_a}(r) \\
& \hspace*{1.5cm} + \sum_{a<b}\left( \overl{\mcU}_{\vec{k}_a}(t) \overl{\mcU}_{\vec{k}_b}(t) 
-\overl{\mcU}_{\vec{k}_{ab}}(t)\right)\symprod{Y_a}{Y_b}(r) \,. \nonumber
\end{align}

As a second step, given  $z_a(t)$, $z_{bc}(t)$ arbitrary functions of time, we
propose the following inversion procedure
\begin{enumerate}
\item  take  the  functions  $\psi_{N(a,b)}(t)  =  \sqrt{2}  \,  N(a,b)  \cos
(N(a,b)\,t)$,   where    $(a,b)\mapsto   N(a,b)\in\{1,\ldots,N\}$\/   is   an
enumeration of the pairs of integers $(a,b)$, $a<b$.
\item  select the following controls in~\eqref{eq:mechsys:affine},
   \begin{align*}
    u_a(t,q) &= v_a(t,q)  + \frac{1}{\eps} w_a\left(\frac{t}{\eps},t\right) \,
    , \\
    w_a(\tau,t)   &=
    -\sum_{c=1}^{a-1}  \psi_{N(c,a)}  (\tau)  +  \sum_{c=a+1}^{m}  z_{ac}(t)
    \psi_{N(a,c)} (\tau) \, ,
  \end{align*}
  where $v_a(t,q)$ are still to be chosen.
\end{enumerate}

\begin{figure}[tbh!]
    \includegraphics[width=.55\textwidth]{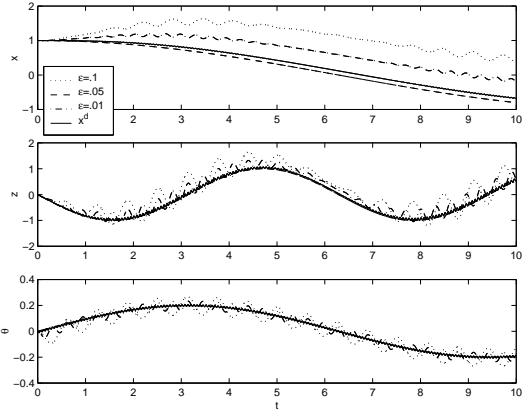}
    \caption{Approximate trajectory tracking for an underactuated PVTOL model
    by means of  oscillatory controls. The curve to be tracked is shown solid,
    and the various oscillating curves correspond to different values of the
    parameter $\eps$}
    \label{fig:pvtol}
\end{figure}

After  computing the  averaged iterated  integrals of  the  oscillatory inputs
$w_a(t/\eps,t)$, equation~\eqref{eq:mechsys:averaged}  for the averaged system
becomes
  \begin{multline*}
    \nabla_{\dot{r}}\dot{r} = Y_0(r)+ k(r)(\dot{r}) +  
     \sum_{a=1}^m v_a(t,r) Y_a(r)  \\ -
    \sum_{a=1}^m \overl{\mcU}_{\vec{k}_{aa}}(t) \symprod{Y_a}{Y_a}(r) 
    + \sum_{a<b} z_{ab}(t)  \symprod{Y_a}{Y_b}(r)\,.
  \end{multline*}
  
  As a third and final step, assume that all the vector fields of the form
  $\symprod{Y_b}{Y_b}$\ belong to $\spn\{Y_a\}$. Let $\alpha_{ab} :
  Q\rightarrow \real$\ be such that $\symprod{Y_a}{Y_a}(q) = \sum_{b}
  \alpha_{ab}(q)Y_b(q)$, $q \in Q$. Select
\[
v_a(t,q) = z_a(t) + \frac{1}{2} \sum_{b=1}^m \alpha_{ba}(q) \left( b-1 +
    \sum_{c=b+1}^m   (z_{bc}(t))^2   \right) \, .
\]
Then, we have
\begin{align*}
  \sum_{a=1}^m  v_a(t,r)Y_a(r) = \sum_{a=1}^m  z_a^d(t) Y_a(r)  + \sum_{a=1}^m
  \overl{\mcU}_{\vec{k}_{aa}}(t)\symprod{Y_a}{Y_a}(r)\,,
\end{align*}
which implies that eq.~\eqref{eq:mechsys:averaged} takes the final form,
  \begin{align*}
    \nabla_{\dot{r}}\dot{r} &= Y_0(r)+ k(r)(\dot{r}) + \sum_{a=1}^m z_a(t) Y_a(r) + 
    \sum_{b<c} z_{bc}(t) \symprod{Y_b}{Y_c}(r)\,,
  \end{align*}
The averaged system  now has more available control  inputs than the original
one.  If the input  distribution $\mcI = \spn \{Y_a\,, \symprod{Y_b}{Y_c}\}$\
is full  rank, then the  latter system is  fully actuated (i.e.,  one control
input is  available for each degree  of freedom).  If  the input distribution
$\mcI$\ contains a  sufficient number of
decoupling vector fields, then  the system is kinematically controllable.  In
both cases, we have reduced the complexity of the motion planning problem.

\begin{figure}[tbh]
    \includegraphics[width=.58\textwidth]{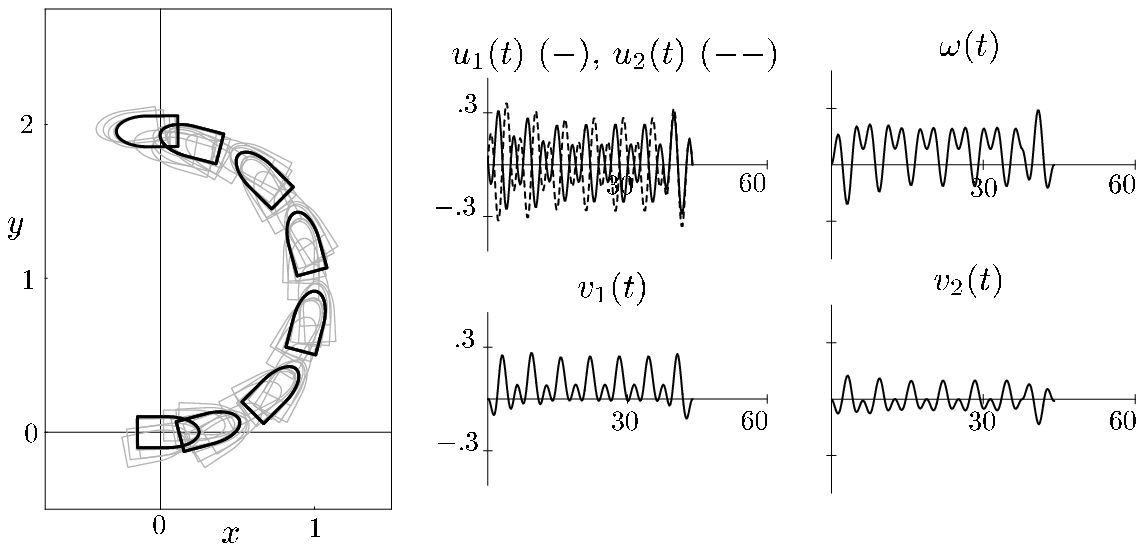}
    \includegraphics[width=.4\textwidth]{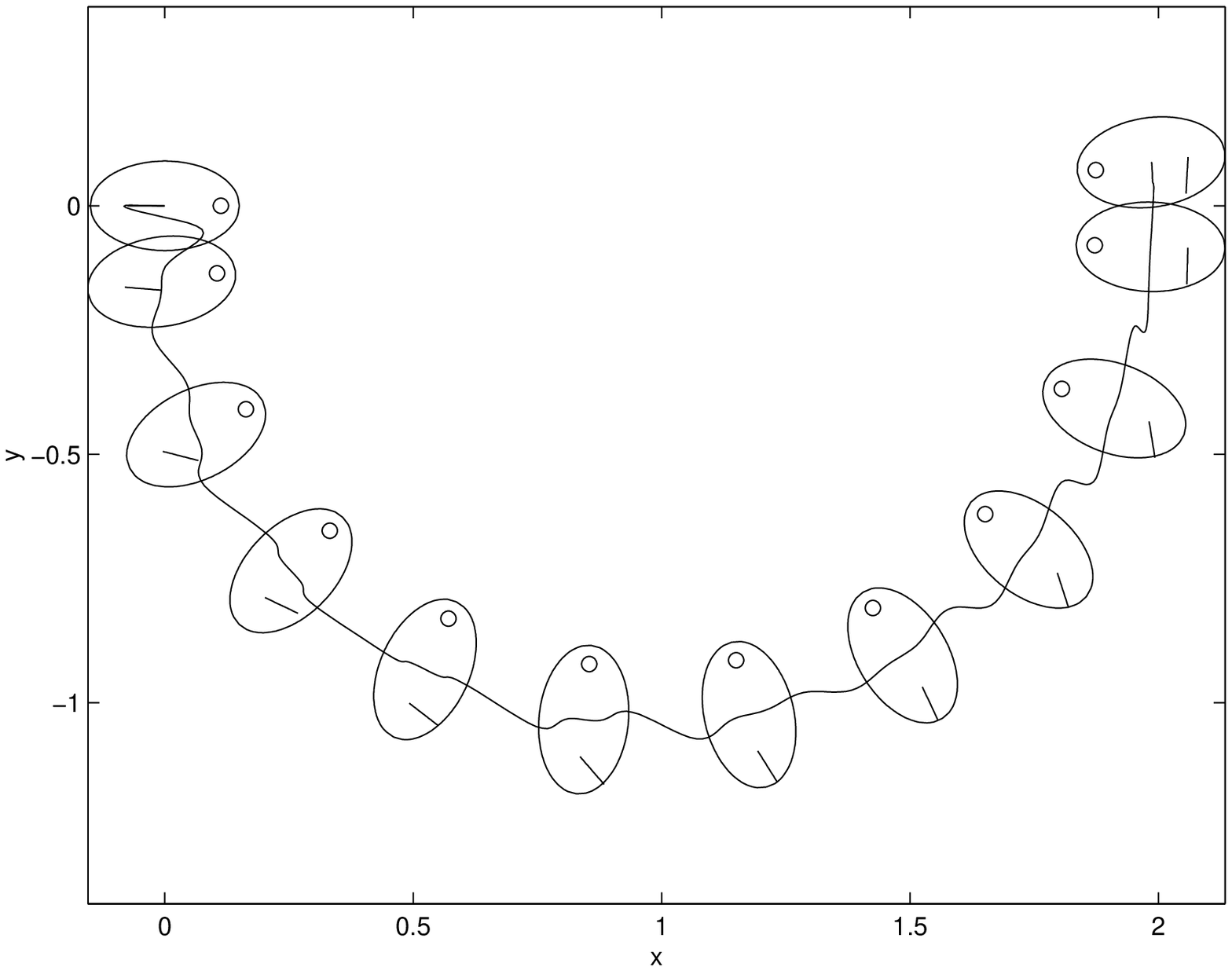}
    \caption{Illustration of the motion planning algorithms via small
    amplitude periodic forcing for a  simple planar body (left) and the blimp
    model (right). The errors in  the final configuration are within the same
    order of magnitude of the input employed}
    \label{fig:small-amplitude}
\end{figure}

\begin{remark}[Small amplitude algorithms based on series expansions]
  {\rm A related approach to motion planning relies on small amplitude periodic
  forcing; see~\cite{FB-NEL-ADL:00,SM-JC:02}. The planning problem is solved
  by approximately inverting the series expansion describing the evolution of
  the control system (cf.  Section~\ref{sec:series-expansion}).  This
  inversion procedure is very similar to the one presented above. Based on
  it, one can establish two simple primitives of motion to change and
  maintain velocity, while keeping track of the changes in the configuration.
  These primitives can then be used as the building blocks to design
  high-level motion algorithms that solve the point-to-point reconfiguration
  problem, the static interpolation problem and the local exponential
  stabilization problem.  Fig.~\ref{fig:small-amplitude} shows two examples
  of the execution of these algorithms.
  }
\end{remark}

\section{Open problems and possible approaches}

Immediate  open questions  arising from  the above-presented  results  are the
following:

\paragraph{Kinematic modeling and control.}
The current limitations are as follows: the design problem is now reduced to
planning for a kinematic system with the additional constraint of
zero-velocity transitions between feasible motions.  This additional
constraint leads to poor performance when coupled with current randomized
planners~\cite{DH-JCL-RM:99a,LEK-PS-JCL-MHO:96,SML-JJK:00,SML-JJK:01} that
switch frequently between the available motions.  The zero-velocity switches
also create problems for trajectory tracking controllers based on
linearization, since the system loses linear controllability at
zero-velocity.  Finally, there is no notion of time-optimality for these
kinematic motions and there is no way of dealing with systems where
oscillatory inputs are needed for locomotion (see below for a
discussion on this point).  Motivated by this analysis, we identify the
following open issues:
\begin{enumerate}
\item Develop a catalog of kinematically controllable systems, including
  planar manipulators with revolute as well as prismatic joints, parallel
  manipulator, manipulators in three dimensional space and in aerospace and
  underwater environments (accounting for the different dynamics in such
  settings).  Some preliminary work in this direction can be found in
  \cite{FB-ADL-KML:02a}.  Analyze and classify the singularities that these
  vector fields possess as a prerequisite step for planning purposes.
  
\item A (left) group action is a map $\psi: G \times Q \rightarrow Q$ such
  that $\psi (e,q)=q$, for all $q \in Q$, where $e$ denotes the identity
  element in $G$, and $\psi(g,\psi(h,q))=\psi(gh,q)$, for all $g,h \in G$, $q
  \in Q$. Usually $G \subset SE(n)$, and then the action describes a rigid
  displacement of some components of the robot.  An interesting problem would
  be to identify conditions under which decoupling vector fields can be found
  which are invariant under such group actions.  When this is the case,
  motion plans can be designed exploiting established ``inverse kinematics''
  methods; see~\cite[Chapter~3]{RMM-ZXL-SSS:94}.  This simplification
  eliminates the need for any numerical procedure if the robot moves in an
  un-obstructed environment, or further reduces the dimensionality and
  complexity of the resulting search problem in complex environments.
  
\item To tackle the difficulties inherent with zero-velocity transitions, it
  would be appropriate to develop randomized planners which require as few
  switches between decoupling vector fields as possible, and to develop
  trajectory tracking controllers for these systems able to adequately
  perform through the singularities.
  
\item Another interesting idea  would consist of switching between decoupling
vector fields  without stopping. In some  sense, this is also  related to the
problem    of   developing    transitions   between    relative   equilibria.
\emph{Relative equilibria} are ``steady trajectories'' that the system admits
as feasible solutions.   This family of trajectories is  of great interest in
theory and  applications as they  provide a rich  family of motions  with the
simplifying  property  of  having  constant  body-fixed  velocity.   Relative
equilibria for systems in  three dimensional Euclidean space include straight
lines, circles, and helices. Despite  partial results, no method is currently
available  to design provably  stable switching  maneuvers from  one relative
equilibrium  to another  (or  from  one decoupling  vector  field to  another
without stopping).  A necessary preliminary  step toward this objective is to
analyze the controllability properties  of underactuated systems moving along
a relative equilibrium or along a decoupling vector field.
 
\end{enumerate}

\paragraph{Small-amplitude and high-frequency controls.}
The current limitations are as follows. The implementation of the small
amplitude approach requires the computation and manipulation of high order
tensors, and the approach has a limited region of convergence.  The
implementation of the oscillatory control approach presents difficulties in
most physical settings because of the required high frequency, high amplitude
inputs.  Motivated by this analysis, we think that the following are
interesting issues to explore:

\begin{enumerate}
 
\item For the small amplitude controls formulation, open questions include
  (a) investigate tight estimates for the region of validity of the
  truncations (simulation studies suggest that there are better bounds than
  the conservative ones currently available), (b) design base functions
  optimal with regards to region of convergence and appropriate cost
  criteria, (c) design inversion algorithms for systems that are not linearly
  controllable.  The latter setting is equivalent to a non-definite quadratic
  programming problem, i.e., to the problem of finding sufficient conditions
  for a vector-valued quadratic form to be surjective
  (see~\cite{FB-JC-ADL-SM:02d} for a discussion on this
  subject).
  
  
\item For the oscillatory controls formulation, standing problems are (a)
  investigate the use of high-frequency bounded amplitude controls, (b)
  characterize approximate kinematic controllability and differential
  flatness via oscillations, (c) investigate physical settings in which
  oscillatory controls are natural control means, e.g.,
  micro-electromechanical robots, (d) investigate extensions of this
  coordinate-free perturbation theory to discrete-time nonlinear systems, and
  to distributed parameter systems and partial differential equations.
  
\item An ambitious program would consist of developing schemes that combine
  the proposed analytic methods with iterative numerical algorithms.  One
  approach is via homotopy and level set methods~\cite{ELA-KG:90,JAS:96} as
  schemes that overcome the limitations induced by the small parameter (small
  convergence region or high amplitude high frequency).  A second direction
  is to use the planner based on small amplitude controls as a local planner
  inside a global search algorithm based on randomization;
  see~\cite{TK-FB:01h} for some preliminary results on local/global planners.
  
\end{enumerate}


\section*{Acknowledgments}
This research was partially funded by NSF grants CMS-0100162 and IIS-0118146.
The authors would like to thank Andrew Lewis, Kevin Lynch, Milo\v s \v
Zefran, Todd Cerven, and Timur Karatas.


\begin{thebibliography}{10}
\addcontentsline{toc}{section}{References}

\bibitem{ELA-KG:90}
Allgower~EL, Georg~K (1990) 
\newblock {\em Numerical Continuation Methods: An Introduction}.
\newblock Springer Verlag, New York.

\bibitem{JBo:00}
Borenstein~J (2000)
\newblock The {OmniMate}: a guidewire- and beacon-free {AGV} for highly
  reconfigurable applications.
\newblock {\em Int J Production Res} 38(9):1993--2010.

\bibitem{FB:99a}
Bullo~F (2001) 
\newblock Series expansions for the evolution of mechanical control systems.
\newblock {\em SIAM J. Control Optim} 40(1):166--190.

\bibitem{FB:99b}
Bullo~F (2002)
\newblock Averaging and vibrational control of mechanical systems.
\newblock {\em SIAM J. Control Optim} 41(2):542--562.

\bibitem{FB-JC-ADL-SM:02d}
Bullo~F, Cort{\'e}s~J, Lewis~AD, Mart{\'\i}nez~S (2002)
\newblock Vector-valued quadratic forms in control theory.  
\newblock In: {\em Mathematical Theory of Networks and Systems (MTNS)}, Notre 
Dame, Indiana, August 2002.  Workshop on Open Problems in Mathematical Systems
and Control Theory.

\bibitem{FB-NEL-ADL:00}
Bullo~F, Leonard~NE, Lewis~AD (2000)
\newblock Controllability and motion algorithms for underactuated {L}agrangian
  systems on {L}ie groups.
\newblock {\em IEEE Trans Automat Control} 45(8):1437--1454.

\bibitem{FB-ADL:00d}
Bullo~F, Lewis~AD (2000)
\newblock On the homogeneity of the affine connection model for mechanical
  control systems.
\newblock In: {\em {IEEE} Int Conf on Decision and Control}, Sydney,
Australia, December 2000, pp 1260--1265.

\bibitem{FB-ADL-KML:02a}
Bullo~F, Lewis~AD, Lynch~KM (2002)
\newblock Controllable kinematic reductions for mechanical systems: concepts,
  computational tools, and examples.
\newblock In: {\em Mathematical Theory of Networks and Systems (MTNS)}, Notre
  Dame, Indiana, August 2002.

\bibitem{FB-KML:01a}
Bullo~F, Lynch~KM (2001)
\newblock Kinematic controllability for decoupled trajectory planning in
  underactuated mechanical systems.
\newblock {\em IEEE Trans Robotics Automat} 17(4):402--412.

\bibitem{FB-MZ:01d}
Bullo~F, {\v Z}efran~M (2002)
\newblock On mechanical control systems with nonholonomic constraints and
  symmetries.
\newblock {\em Systems \& Control Lett} 45(2):133--143.

\bibitem{MC-FSL-PEC:95}
Camari{\~n}a~M, Silva Leite~F, Crouch~PE (1995)
\newblock Splines of class $C^k$ on non-{E}uclidean spaces.
\newblock {\em IMA J Math Control \& Information} 12:399--410.

\bibitem{GC-GB-BDAN:96}
Campion~G, Bastin~G, D'Andrea-Novel~B (1996)
\newblock Structural properties and classification of kinematic and dynamic
  models of wheeled mobile robots.
\newblock {\em IEEE Trans Robotics Automat} 12(1):47--62.

\bibitem{GC-JWB:94}
Chirikjian~G, Burdick~JW (1994)
\newblock Kinematics of hyperredundant locomotion.
\newblock {\em IEEE Trans Robotics Automat} 11(6):781--793.

\bibitem{JC-SM-FB:01m}
Cort{\'e}s~J, Mart{\'\i}nez~S, Bullo~F (2002)
\newblock On nonlinear controllability and series expansions for {L}agrangian
  systems with dissipative forces.
\newblock {\em IEEE Trans Automat Control} 47(8):1396--1401.

\bibitem{GE-SH:99}
Endo~G, Hirose~S (1999)
\newblock Study on roller-walker (system integration and basic experiments).
\newblock In: {\em {IEEE} Int Conf on Robotics and Automation}, Detroit,
Michigan, May 1999, pp 2032--2037.

\bibitem{SH:93}
Hirose~S (1993)
\newblock {\em Biologically inspired robots: snake-like locomotors and
  manipulators}.
\newblock Oxford University Press, Oxford.

\bibitem{SH:00}
Hirose~S (2000)
\newblock Variable constraint mechanism and its application for design of
  mobile robots.
\newblock {\em Int J Robotics Res} 19(11):1126--1138.

\bibitem{JKH-MHR:90}
Hodgins~JK, Raibert~MH (1990)
\newblock Biped gymnastics.
\newblock {\em Int J Robotics Res} 9(2):115--132.

\bibitem{RH-OK:00}
Holmberg~R, Khatib~O (2000)
\newblock Development and control of a holonomic mobile robot for mobile
  manipulation tasks.
\newblock {\em Int J Robotics Res} 19(11):1066--1074.

\bibitem{DH-JCL-RM:99a}
Hsu~D, Latombe~JC, Motwani~R (1999)
\newblock Path planning in expansive configuration spaces.
\newblock {\em Int J Comp Geom Appl} 9(4):495--512.

\bibitem{JCJ-SK-JA:95}
Jalbert~JC, Kashin~S, Ayers~J (1995)
\newblock Design considerations and experiments of a biologically based
  undulatory lamprey {AUV}.
\newblock In: {\em 9th Int Symposium on Unmanned Untethered
  Submersible Technology}, Durham, New Hampshire.

\bibitem{TK-FB:01h}
Karatas~T, Bullo~F (2001)
\newblock Randomized searches and nonlinear programming in trajectory planning.
\newblock In: {\em {IEEE} Conf on Decision and Control}, Orlando, Florida,
December 2001, pp 5032--5037.

\bibitem{NK-TI:98}
Kato~N, Inaba~T (1998)
\newblock Guidance and control of fish robot with apparatus of pectoral fin
  motion.
\newblock In: {\em {IEEE} Int Conf on Robotics and Automation}, Leuven,
Belgium, May 1998, pp 446--451.

\bibitem{LEK-PS-JCL-MHO:96}
Kavraki~LE, {\v S}vestka~P, Latombe~JC, Overmars~MH (1996)
\newblock Probabilistic roadmaps for path planning in high-dimensional space.
\newblock {\em IEEE Trans Robotics Automat} 12(4):566--580.

\bibitem{SDK-RJM-CTA-RMM-JWB:98}
Kelly~SD, Mason~RJ, Anhalt~CT, Murray~RM, Burdick~JW (1998)
\newblock Modelling and experimental investigation of carangiform locomotion
  for control.
\newblock In: {\em {IEEE} {A}merican {C}ontrol {C}onference}, Philadelphia,
PA, 1998, pp 1271--1276.

\bibitem{PSK-DPT:94}
Krishnaprasad~PS, Tsakiris~DP (1994)
\newblock {G}-snakes: Nonholonomic kinetic chains on {Lie} groups.
\newblock In: {\em {IEEE} Int Conf on Decision and Control}, Lake
  Buena Vista, Florida, December 1994, pp 2955--2960.

\bibitem{PSK-DPT:01}
Krishnaprasad~PS, Tsakiris~DP (2001)
\newblock Oscillations, ${SE}(2)$-snakes and motion control: a study of the
  roller racer.
\newblock {\em Dynam Stab Systems} 16(4):347--397.

\bibitem{SML-JJK:00}
LaValle~SM, Kuffner~JJ
\newblock Rapidly-exploring random trees: Progress and prospects.
\newblock In: {\em Workshop on Algorithmic Foundations of Robotics},
Dartmouth, New Hampshire, March 2000, pp 293--308.

\bibitem{SML-JJK:01}
LaValle~SM, Kuffner~JJ. (2001)
\newblock Randomized kinodynamic planning.
\newblock {\em Int J Robotics Res} 20(5):378--400.

\bibitem{ADL:99a}
Lewis~AD (2000)
\newblock The geometry of the maximum principle for affine connection control
  systems.
\newblock Preprint available at:
\texttt{http:/\!/penelope.mast.queensu.ca/\~{}andrew}. 

\bibitem{ADL:97a}
Lewis~AD (2000)
\newblock Simple mechanical control systems with constraints.
\newblock {\em IEEE Trans Automat Control} 45(8):1420--1436.

\bibitem{ADL-RMM:95c}
Lewis~AD, Murray~RM (1997)
\newblock Configuration controllability of simple mechanical control systems.
\newblock {\em SIAM J Control Optim} 35(3):766--790.

\bibitem{ADL-JPO-RMM-JWB:94}
Lewis~AD, Ostrowski~JP, Murray~RM, Burdick~JW (1994)
\newblock Nonholonomic mechanics and locomotion: the snakeboard example.
\newblock In: {\em {IEEE} Int Conf on Robotics and Automation}, San Diego,
California, May 1994, pp 2391--2400.

\bibitem{SM-JC:02}
Mart{\'\i}nez~S, Cort\'es~J (2002)
\newblock Motion control algorithms for simple mechanical systems with
  symmetry.
\newblock {\em Acta Applicandae Mathematicae}, to appear.

\bibitem{SM-JC-FB:01e}
Mart{\'\i}nez~S, Cort{\'e}s~J, Bullo~F (2001)
\newblock Analysis and design of oscillatory controls systems.
\newblock {\em IEEE Trans Automat Control}, submitted. Available
electronically at \texttt{http:/\!/motion.csl.uiuc.edu}.

\bibitem{TM:90}
McGeer~T (1990)
\newblock Passive dynamic walking.
\newblock {\em Int J Robotics Res} 9(2):62--82.

\bibitem{KAM-JPO:99a}
Mclsaac~KA, Ostrowski~JP (1999)
\newblock A geometric approach to anguilliform locomotion: modeling of an
  underwater eel robot.
\newblock In: {\em {IEEE} Int Conf on Robotics and Automation}, Detroit,
Michigan, May 1999, pages 2843--2848.

\bibitem{RMM-ZXL-SSS:94}
Murray~RM, Li~ZX, Sastry~SS (1994)
\newblock {\em A Mathematical Introduction to Robotic Manipulation}.
\newblock CRC Press, Boca Raton, Florida.

\bibitem{SO-JA-RS:00}
Ostrovskaya~S, Angeles~J, Spiteri~R (2000)
\newblock Dynamics of a mobile robot with three ball-wheels.
\newblock {\em Int J Robotics Res} 19(4):383--393.

\bibitem{JPO:00a}
Ostrowski~JP (2000)
\newblock Steering for a class of dynamic nonholonomic systems.
\newblock {\em IEEE Trans Automat Control} 45(8):1492--1497.

\bibitem{JPO-JWB:98}
Ostrowski~JP, Burdick~JW (1998)
\newblock The geometric mechanics of undulatory robotic locomotion.
\newblock {\em Int J Robotics Res} 17(7):683--701.

\bibitem{FGP-SMK:94}
Pin~FG, Killough~SM (1994)
\newblock A new family of omnidirectional and holonomic wheeled platforms for
  mobile robots.
\newblock {\em IEEE Trans Robotics Automat} 10(4):480--489.

\bibitem{MHR:86}
Raibert~MH (1986)
\newblock {\em Legged Robots that Balance}.
\newblock MIT Press, Cambridge, Massachusetts.

\bibitem{MR-RMM:96b}
Rathinam~M, Murray~RM (1998)
\newblock Configuration flatness of {L}agrangian systems underactuated by one
  control.
\newblock {\em SIAM J Control Optim} 36(1):164--179.

\bibitem{SS-JA-JD:95}
Saha~S, Angeles~J, Darcovich~J (1995)
\newblock The design of kinematically isotropic rolling robots with
  omnidirectional wheels.
\newblock {\em Mechanism and Machine Theory} 30(8):1127--1137.

\bibitem{JAS:96}
Sethian~JA (1996)
\newblock {\em Level set methods: {E}volving interfaces in geometry, fluid
  mechanics, computer vision, and materials science}.
\newblock Cambridge University Press, New York.

\bibitem{MW-HA:97}
West~M, Asada~H (1997)
\newblock Design of ball wheel mechanisms for omnidirectional vehicles with
  full mobility and invariant kinematics.
\newblock {\em ASME J Mech Design} 119(2):153--161.

\end{thebibliography}

\end{document}